\documentclass[a4 paper]{article}
\usepackage{amsthm, a4wide}
\usepackage{amsfonts}
\usepackage{amsmath}
\usepackage{amssymb}
\usepackage[english]{babel}
\newtheorem{defi}{Definition}
\newtheorem{teo}[defi]{Theorem}
\newtheorem{rem}[defi]{Remark}

\newtheorem{lemma}[defi]{Lemma}

\newcommand{\esup}{\text{ess}\sup_{\omega\in\Omega}}
\newcommand{\eps}{\varepsilon}
\newcommand{\ox}{\overline{X}}
\newcommand{\oa}{\overline{\alpha}}
\newcommand{\R}{I\!\!R}
\begin{document}
\title{A converse Lyapunov theorem for
almost sure stabilizability
\thanks{This research was partially supported by M.I.U.R.,
project ``Viscosity, metric, and control theoretic methods for nonlinear
partial differential equations'', and by GNAMPA-INDAM, project ``Partial
differential equations and control theory''.}}\author{ Annalisa Cesaroni\\
Dipartimento di Matematica P.  e A.\\ Universit\`a di Padova\\
via Belzoni 7, 35131 Padova, Italy\\
acesar@math.unipd.it\\ \\
}

\date{}
\maketitle
\begin{abstract}
We prove a converse Lyapunov theorem for almost sure
stabilizability and almost sure asymptotic stabilizability of
controlled diffusions: given a stochastic system a.s.  stochastic
open loop   stabilizable at the origin, we construct a lower
semicontinuous positive definite function whose level sets form a
local basis of viable neighborhoods of the equilibrium. This
result provides, with the direct Lyapunov theorems proved in a
companion paper, a complete Lyapunov-like characterization of the
a.s. stabilizability.


%\smallskip \noindent{\bf AMS subject classification}. 34D20, 93E15, 49L25.
\end{abstract}

\section{Introduction}
In this paper we provide a Lyapunov characterization of almost
sure  stochastic open loop   stability  at an equilibrium of
controlled diffusion processes in $\R^{N}$
\[(CSDE)\left\{\begin{array}{l}
dX_t=f(X_t,\alpha_t)dt
+\sigma(X_t,\alpha_t)dB_t,\;\;\alpha_t\in A,\;\; t>0,\\ X_0=x,
\end{array}
\right. \] This notion has been introduced in a companion paper
\cite{bc2} by Bardi and the author (see also ~\cite{bc1}): we say
that $(CSDE)$ is {\em a.s. (open loop) stabilizable}  if  for any
$\eta>0$ there exists $\delta>0$ such that, for any $x$ with
$|x|\leq \delta$, there exists $\alpha$ such that the
corresponding process satisfies $|X_t|\leq \eta$ for all $t\geq 0$
almost surely. If, in addition, the trajectory is asymptotically
approaching a.s. the equilibrium, we say the system is {\em a.s.
(open loop) asymptotically stabilizable}. The definitions imply in
particular that these properties are never verified by
nondegenerate processes. This stochastic stability describes a
behaviour very similar to a stable deterministic system and is
stronger than pathwise stability and stability in probability (see
\cite{has, kus, maolib}). We characterize it by means of
appropriate control {\em Lyapunov functions}. These functions have
been introduced in \cite{bc1} and are lower semicontinuous (LSC),
continuous at the equilibrium, positive definite, proper. Moreover
they satisfy the following
 {\em infinitesimal decrease condition}:  \begin{equation}
\label{cseqintro}\sup_{ \sigma(x,\alpha)^T DV(x) = 0}\!\!\!
\left\{-DV(x)\cdot f(x,\alpha)-trace\left[a(x,\alpha)
D^2V(x)\right] \right\} \geq l(x),
\end{equation}
 where $a:=\sigma\sigma^T/2$, $l\equiv 0$ for
Lyapunov functions and positive definite for strict Lyapunov
functions. This is not a standard Hamilton-Jacobi-Bellman
inequality, because the constraint on the controls depends on $V$:
we are allowing diffusion only in the directions tangential to the
sublevel sets of $V$. If we eliminate this constraint, the
differential inequality which is left is the infinitesimal
decrease condition on  Lyapunov functions for the stability in
probability (see \cite{thesis},\cite{ces2}). We prove that $V$
satisfies (\ref{cseqintro}) if and only if it satisfies the
following {\em monotonicity condition:}
\begin{equation}\label{mono}\forall x\ \exists \alpha\ :\quad\sup_{t\geq
0}\esup\left(V(X_t)+\int_0^t l(X_s)ds\right)\leq
V(x)\end{equation} where the essential supremum is intended  with
respect to the probability measure $\mathbf{P}_x$. This means in
particular that the process $V(X_t)$ is a positive {\em
supermaxingale} according to the definition given in \cite{bcj}:
this is the natural counterpart of the requirement on the process
$V(X_t)$ to be a positive supermartingale in the context of
stability in probability. The monotonicity condition says
 that the sublevel sets  $K_\mu:=\{x\
\ V(x)\leq \mu\}$ are {\em viable (or weakly invariant)} with
respect to (CSDE) in the sense that $\forall \ x\in K_\mu\ \exists
\alpha$ such that $X_t\in K_\mu$ forever almost surely. One of the
main tool used in this  paper  is the geometric Nagumo-type
characterization of viability proved recently by Bardi and Jensen
in \cite{bj} (see also \cite{bg}, \cite{ad3} and the references
therein for earlier related results).

In \cite{bc2}, Bardi and the author show that the existence of a
Lyapunov function (respectively, of a strict Lyapunov function)
implies the a.s. stabilizability (respectively, the a.s.
asymptotic stabilizability) of the system to the equilibrium. As a
simple example of application of this theory, we consider a radial
function $V(x)=v(|x|)$, for some real smooth function $v$ with
$v'(r)>0$ for $r>0$. The system (CSDE) admits $V$ as Lyapunov
function if \[ \forall x \ \exists \alpha:\ \sigma(x,\alpha)\cdot
x = 0 \quad f(x,\alpha)\cdot x+trace\ a(x,\alpha)\leq 0.\]
Therefore the following  conditions are sufficient for the a.s.
stabilizability: the radial component of the diffusion is null and
its rotational component, which still plays a destabilizing role
since $trace\, a(x,\alpha)\geq 0$, must be compensated by a
negative radial component of $f$.

In this paper we prove that  the existence of a Lyapunov function
is also a necessary condition for a.s. stabilizability (and also
for a.s. Lagrange stabilizability in the global case). A Lyapunov
function for the system can be defined as
\[V(x):=\inf\{r\ |\ \exists \oa \mbox{ admissible control such that }
|\ox_t|\leq r \mbox{ almost surely } \forall t\geq 0\},\] or
equivalently as
\[V(x):=\inf_{\alpha\in \mathcal{A}_x}\sup_{t\geq 0}\mbox{ess}\sup_{\omega\in \Omega}
|X_t^\alpha|.\] We prove that this function is LSC, continuous at
the origin, positive definite, proper,
 and satisfies the infinitesimal decrease condition (\ref{cseqintro}) with $l\equiv 0$.
For an a.s. asymptotic stabilizable systems in a bounded set
$\mathcal{O}$, we  build a positive definite Lipschitz continuous
function $l$, related to  the rate of decrease of the stable
trajectories to the equilibrium, by the formula
\[V(x):=\inf_{\alpha\in\mathcal{A}_x}\esup\
\int_0^{+\infty}l(X_t^\alpha)dt.\]We show that $V$ is finite, LSC,
continuous at the origin, positive definite, proper,
 and satisfies the infinitesimal decrease condition
 (\ref{cseqintro}).

In both cases the Lyapunov functions are  {\em worst-case value
function} of an appropriate stochastic optimal control problem: we
minimize the worst possible cost over all possible paths. This is
quite natural since these functions characterize a very strong
stability notion. The link between worst-case value functions and
viscosity solutions to geometric second order partial differential
equations has been recently treated  by Soner and Touzi in
~\cite{st2} (see also \cite{bcj}). They considered stochastic
target problems where the controller tries to steer almost surely
a controlled process into a given target by judicial choices of
controls. The interest in this kind of stochastic control problems
in the {\em almost sure } setting comes from the relationship with
mean curvature type geometric flows and from the applications to
the super-replication problems in financial mathematics. Moreover
the a.s. stability of control systems affected by disturbances
modelled as  $M$-dimensional white noise  is related to the
so-called {\em worst-case stability (or robust stability)} of
deterministic control systems affected by disturbances modelled as
(deterministic) $L^\infty$ functions with values in $\R^M$ (see
\cite{fr}). The Lyapunov characterization of these two stability
properties seems to be an useful tool to give a precise proof of
this relationship (see \cite{thesis}), while a direct proof  based
on the estimates among the trajectories of the two systems should
be rather hard. An approach of this type has been used recently by
Da Prato and Frankowska in \cite{af} to prove the equivalence
between the invariance with respect to a controlled stochastic
system and the invariance with respect to a deterministic system
with two (non competitive) controls.

We conclude with some additional references on converse Lyapunov
theorems. For  controlled deterministic systems, there are
theorems characterizing the  stochastic open loop stabilizability
by means of LSC appropriate Lyapunov functions (see \cite{br}).
Soravia in \cite{sor1}  showed  that the stability at an
equilibrium is equivalent to the continuity at such point of the
value function $V(x)=\inf_{\alpha\in\mathcal{A}_x}\sup_{t\geq 0}
U(X_t)$ where the level sets of $U$ form a local basis of
neighborhoods of the equilibrium. For asymptotically controllable
systems, Sontag and Sussmann (\cite{son0}, \cite{sosu}) provided a
characterization of asymptotic controllability by means of
continuous Lyapunov functions such as
$V(x)=\inf_\alpha\int_0^{+\infty}l(X_t)dt$ where $l$ is an
appropriate positive definite function. Recently Rifford
(\cite{r}) proved a converse Lyapunov theorem in the framework of
Lipschitz continuous functions which are semiconcave outside the
equilibrium. In the stochastic setting, Has'minskii (\cite{has},
\cite{kus}) obtained a converse theorem for  stability in
probability
 of  uncontrolled diffusion processes,
strictly nondegenerate outside the equilibrium, by means of
$\mathcal{C}^2$ Lyapunov functions, using the Maximum Principle
and the properties of solutions of uniformly elliptic equation.
Kushner proved in \cite{kus1} a characterization of asymptotic
uniform stochastic stability by means of continuous Lyapunov
functions (here, however, the infinitesimal decrease condition is
given in terms of the weak generator of the process). In the
forthcoming paper \cite{ces2} (see also \cite{thesis}) the author
extends the direct Lyapunov method  by Has'minskii and Kushner to
the study of stochastic open loop  stabilizability in probability
in terms of merely semicontinuous Lyapunov functions which satisfy
in the viscosity sense an appropriate infinitesimal decrease
condition and provides also in this setting converse Lyapunov
theorems.

The paper is organized as follows. In Section 2 we  give the
definition of  stochastic open loop   a.s. stabilizability, in
Section 3 we introduce the appropriate concept of Lyapunov
function for the study of such stability. Section 4 is devoted to
the viability properties of sublevel sets of Lyapunov functions.
Section 5 contains the main results: the converse Lyapunov
theorems. Finally in Section 6 we give the extension to general
equilibrium sets.

\section{Almost sure  Lyapunov stabilizability}

We consider a controlled  Ito stochastic differential equation:
\[(CSDE)\left\{\begin{array}{l}
dX_t=f(X_t,\alpha_t)dt
+\sigma(X_t,\alpha_t)dB_t,\;\; t>0,\\ X_0=x.
\end{array}
\right. \]
where $B_t$ is an $M$-dimensional Brownian motion. Throughout the paper we
assume that
$f, \sigma$ are continuous functions defined in
$\R^N\times A$, where $A$ is a compact metric
space, which take values, respectively, in $\R^N$ and in the space of
$N\times M$ matrices, and satisfying for all $x,y\in\R^N$ and all $\alpha\in A$
\begin{equation}\label{condition1}
\mbox{$|f(x,\alpha)-f(y,\alpha)|+
\Vert\sigma(x,\alpha)-\sigma(y,\alpha)\Vert\leq C|x-y|$,
}
\end{equation}
\noindent We define
$a(x,\alpha):=\frac{1}{2}\sigma(x,\alpha)\sigma(x,\alpha)^T$ and
assume
\begin{equation}\label{convex}
\left\{(a(x,\alpha),f(x,\alpha)) \; : \; \alpha\in A \right\}\quad
\text{is convex for all } x\in \R^N .
\end{equation}
The class of admissible controls is the class of {\em strict
controls}, as defined in ~\cite[Definition 2.2]{hl}: they are $A$
valued, progressively measurable processes $\alpha_t$ such that
there exists a solution $X_t^\alpha$ to $(CSDE)$. $\mathcal{A}_x$
denotes  the class of admissible control for a given initial datum
$x$, with $\alpha_\cdot$ its generic element (although it is not a
standard function $\R\to A$), and with $X_\cdot$ the corresponding
solution of $(CSDE)$. We recall also a theorem on the existence of
optimal control for stochastic control problems.
\begin{teo}[Theorem 4.7 and Corollary 4.8 \cite{hl}] \label{optimalcontrol}
Under the convexity assumption (\ref{convex}), for every initial
data $x\in\R^N $ there exists an admissible control realizing the
minimum in  the control problem
$\inf_{\alpha}\mathbf{E}J(x,\alpha)$ where the cost functional
$J(x,\alpha)$ satisfies standard regularity assumptions.
\end{teo}
We state  now the definition of almost sure  stochastic open loop
stabilizability, which has been introduced and studied in the
paper ~\cite{bc2} (see also ~\cite{bc1} for the uncontrolled
case). We introduce the classes of comparison functions.
\begin{defi}[comparison functions]\upshape $\mathcal{K}$
denotes the class of real continuous functions $\gamma$ strictly
increasing and such that $\gamma(0)=0$; $\mathcal{K}_{\infty}$
contains the functions $\gamma\in\mathcal{K}$ such that
$\lim_{r\rightarrow +\infty}\gamma(r)=+\infty$. Finally
$\mathcal{KL}$ denotes the class of continuous functions
$\beta:\R\times \R\to\R$ which are strictly increasing in the
first variable, strictly decreasing in the second variable, which
satisfy $\beta(0,t)=0$ for  $t\geq 0$, $\lim_{t\rightarrow
+\infty}\beta(r,t)=0$ for $r\geq 0.$
\end{defi}
\begin{defi}[a.s. stabilizability]\label{def:as.stab}
The system $(CSDE)$ is \emph{almost surely (open-loop Lyapunov) stabilizable}
%or \emph{a.s.  asymptotically controllable}
at the origin if there exists $\gamma\in\mathcal{K}$ and
$\delta_o>0$ such that for any starting point $x$ with $|x|\leq
\delta_o$
\begin{equation}
\label{cappa}
 %\exists \gamma\in\mathcal{K} : \; \forall x\in\R^N \;
\exists \, \overline{\alpha}_{\cdot}\in{\cal A}_x \,: \quad
|\overline{X}_t|\leq \gamma(|x|)\quad \forall t\geq 0 \,\text{
a.s.}
\end{equation}
If $\gamma$ can be chosen in $\mathcal{K}_{\infty}$ and the
estimate (\ref{cappa}) holds in the whole space $\R^N$, the system
is also \emph{almost surely (open-loop) Lagrange stabilizable}, or
it has the property of \emph{uniform boundedness of trajectories}.
\end{defi}
\begin{rem}\upshape
We could define the a.s. stabilizability equivalently as follows:

\noindent {\it the system is a.s. (open loop) stabilizable at the
origin if  for every  $\eps>0$ there exists $\delta>0$  such that
for $|x|\leq\delta$  there exists an admissible control function
$\overline{\alpha}_{\cdot}\in{\cal A}_x$ whose corresponding
trajectory $\overline{X}_{\cdot}$ verifies $|\overline{X}_t|\leq
\eps$ for all $t\geq 0$ almost surely.}
\end{rem}
\begin{rem}\upshape
A necessary condition for the a.s. stabilizability at the  origin
is that the origin  is a {\em controlled equilibrium} of $(CSDE)$,
i.e.,
%for some $\oa\in A$ .
\begin{equation} \label{equilibrium}
\exists \,\oa\in A \,:\; f(0,\oa)=0, \; \sigma(0,\oa)=0.
\end{equation}
\end{rem}
\begin{defi}[a.s. asymptotic stabilizability]\label{asstabasinto}
The system $(CSDE)$ is \emph{almost surely (open loop) locally
asymptotically stabilizable} (or \emph{a. s. locally
asymptotically controllable}) at the origin if  there is $\beta\in
\mathcal{KL}$   and $R>0$ such that for any starting point $x$
with $|x|\leq R$ there exists $\overline{\alpha}_{\cdot}\in{\cal
A}_{x}$
\begin{equation}\label{decayest.}|\ox_t|\leq\beta(|x|,t)\quad \forall t\geq 0
\text{ a.s.}\end{equation}\end{defi}

\section{Lyapunov functions for a.s. stabilizability}
In this section we introduce the appropriate concept of Lyapunov
function for the study of the almost sure stochastic stability.\\
We recall the definition of the second order semijet (see
\cite{cil}) of a LSC function $V$ at  $x$
 ${\cal J}^{2,-} V(x) := \{ (p,Y)\in\R^N\times S(N)$ such that for $y\rightarrow
x$,  $V(y)\geq V(x)+p\cdot (y-x)+\frac{1}{2} (y-x)\cdot
Y(y-x)+o(|y-x|^2)\}$.

\begin{defi}[control Lyapunov function]\label{cliap}
Let ${\cal O}\subseteq \R^N$ be an open set containing the origin.
A function $V : {\cal O}\rightarrow[0,+\infty)$ is a \emph{control
Lyapunov function} for the a.s.  stabilizability of $(CSDE)$ if

\noindent (i) $V$ is lower semicontinuous and continuous at $0$;

\noindent(ii) $V$ is   {\em positive definite}, i.e., $V(0)=0$ and
$V(x)>0$ for all $x\neq 0$;

\noindent(iii) $V$ is {\em proper}, i.e., the sublevel sets
$\left\{x| V(x)\leq\mu\right\}$ are bounded $\forall \mu\in
[0,\infty)$;

\noindent (iv) $V$ is a viscosity supersolution in ${\cal
O}\setminus\{0\}$ of the equation:
\begin{equation}\label{cseq1} \sup_{\sigma(x,\alpha)\cdot
DV(x)=0}\!\!\!\!\!\!\left\{-DV(x)\cdot
f(x,\alpha)-trace[a(x,\alpha) D^2V(x)]\right\}\geq
0,\end{equation} in the following sense: for all $x\in{\cal
O}\setminus\{0\}$ and $(p,Y)\in{\cal J}^{2,-}V(x)$ there exists
$\overline{\alpha}\in A$:
\[ \sigma(x,\overline{\alpha})^Tp=0 \;\ \ \  \,\mbox{ and
}\,\ \ \ -p\cdot
f(x,\overline{\alpha})-trace\left[a(x,\overline{\alpha}) Y\right]
\geq 0 .
\]
If there exists a positive definite, Lipschitz continuous $l:{\cal
O}\to \R$ such that
\begin{equation}\label{cseq2} \sup_{\sigma(x,\alpha)\cdot
DV(x)=0}\!\!\!\!\!\!\left\{-DV(x)\cdot
f(x,\alpha)-trace[a(x,\alpha) D^2V(x)]\right\}\geq
l(x)\end{equation}then $V$ is a \emph{strict control Lyapunov
function} for the a.s. stabilizability of $(CSDE)$.
\end{defi}

The inequality (\ref{cseq1})  is not the standard
Hamilton-Jacobi-Bellman inequality arising in stochastic optimal
control. We have an implicit constraint on the controls,
$\sigma(x,\alpha)\cdot DV(x)=0$, i.e. depending on the generalized
subgradients of the solution: we are allowing only controls which
render the diffusion matrix tangential in some generalized sense
to the sublevel sets of  $V$. This implies that the diffusion has
to degenerate in a large set, for some control.

Because of this constraint, though, the nonlinearity
\[F(x,p,X)=\sup\left\{-p\cdot f(x,\alpha)-trace[a(x,\alpha) X]\ |\
a\in A\ \sigma(x,\alpha)\cdot p=0\right\}\]  is  {\bf geometric}
in the sense that it satisfies the following rescaling property
$F(x,\lambda p, \lambda X+\mu p\otimes p)=\lambda F(x,p,X)$ for
every $\lambda>0$ and $\mu\in\R$, where $p\otimes p$  is the
$N\times N$ matrix whose $(i,j)$ entry is $p_ip_j$. This  permits
to prove  the following lemma on the change of unknown (for the
proof see \cite{bc2}, \cite{thesis}).
\begin{lemma}\label{lemma2}
Assume  that  $v$ is a LSC viscosity supersolution of equation
(\ref{cseq1}) in an open set $\mathcal{O}$. Let $\phi$ be a twice
continuously differentiable strictly increasing real map. Then
$w=\phi\circ v$ is still a viscosity supersolution of equation
(\ref{cseq1}) in $\mathcal{O}$.
\end{lemma}
\section{Viability properties of Lyapunov functions}
We study now a viability property of the sublevel sets of
viscosity supersolution of the nonstandard Hamilton-Jacobi-Bellman
inequality (\ref{cseq1}). We recall the definition of almost sure
viability (named also {\em controlled invariance} and {\em weak
invariance}) of an arbitrary closed set for a controlled diffusion
process.
\begin{defi}[viable set] A closed set $K\subset\R^N$ is {\em viable}
or {\em controlled invariant} or {\em weakly invariant} for the
stochastic system $(CSDE)$ if for all initial points $x\in K$
there exists an admissible control $\alpha_.\in{\cal A}_x$ such
that the corresponding trajectory $X_.$ satisfies $X_t\in K$ for
all $t>0$ almost surely.
\end{defi}
This property was studied by Aubin and Da Prato \cite{ad3}
 and, more recently, by Bardi and
Jensen \cite{bj}. The main result of \cite{bj} is the equivalence
between the viability of a closed set $K$ and a Nagumo-type
geometric condition. This geometric condition is  given in terms
of the {\em second order normal cone} to a closed set
$K\subset\R^N$, first introduced in \cite{bg},
\begin{eqnarray*}
{\cal N}^2_K (x)\!\! & := &\!\! \{(p,Y)\in I\!\!R^N \times S(N):
\;\;{\rm for}\;\; y \rightarrow x , \;\;y\in K ,\;\; \\
&    & p\cdot (y-x)+\frac{1}{2} (y-x)\cdot Y(y-x)\geq o(|y-x|^2)
\;\;\}
\end{eqnarray*}
where $S(N)$ is the set of symmetric $N\times N$ matrices.  Note
that, if $\partial K$ is a smooth surface in a neighborhood of
$x$, $p/|p|$ is the interior normal and $Y$ is related to the
second fundamental form of $\partial K$ at $x$, see \cite{bg}.
\begin{teo}[Viability theorem \cite{bj}]
\label{cor.convex}
 Assume conditions (\ref{condition1}) and (\ref{convex}).
Then $K$ is viable for $(CSDE)$ if and only if \[ \forall
x\in\partial K, \; \forall (p,Y)\in {\cal N}^2_K (x),\; \exists
\alpha\in A\;:\;\; f(x,\alpha)\cdot p + trace\left[a(x,\alpha)
Y\right]\geq 0. \] Moreover, for the same $\alpha$ we have that
$\sigma(x,\alpha)\cdot p=0.$
\end{teo}
Using this result we obtain the following characterization, which
is a variant of a result contained in \cite{bc2} (see also
\cite{bc1}).
\begin{lemma}\label{lemma3}Assume conditions (\ref{condition1}) and (\ref{convex}).
Consider an open set $\mathcal{O}\subseteq \R^N$ and a  LSC
function $V:\mathcal{O}\to\R$. If $V$ is a viscosity supersolution
of
\begin{equation}\label{lemma4}
\sup_{\sigma(x,\alpha)\cdot DV=0}\left\{-DV\cdot
f(x,\alpha)-trace[a(x,\alpha)D^2V]\right\}\geq 0\end{equation} in
$\mathcal{O}$, then the sublevel sets $\{V(x)\leq\mu\}$ whose
boundary  is entirely contained in $\mathcal{O}$ are viable with
respect to (CSDE). Viceversa, suppose $\overline{\mu}$ is the
maximal value for which the sublevel set $\{V(x)\leq\mu\}$ has
boundary entirely contained in $\mathcal{O}$. Then the function
$\overline{V}(x):=V(x)\wedge\overline{\mu}$ is a viscosity
supersolution of (\ref{lemma4}) in $\mathcal{O}$.
\end{lemma}
Observe that if $\mathcal{O}=\R^N$, every sublevel set of $V$ is
closed (in particular we can take $\overline{\mu}=+\infty$).
\begin{proof}
For every $\mu\leq \overline{\mu}$, $K_\mu:=\{x \ | \ V(x)\leq
\mu\}$. We define now the LSC function $V_\mu(x):=\mu$ for $x\in
K_\mu$ and $ +\infty $ elsewhere. From the definitions, it is easy
to check that ${\cal J}^{2,-} V_\mu(x)=-{\cal N}^{2}_K(x),\quad
\forall \, x\in
\partial K_\mu,$
so, by the Viability Theorem \ref{cor.convex}, $V_\mu$ is a viscosity supersolution
of (\ref{lemma4}) if and only if $K_\mu$ is viable.

We assume that $V$ is a viscosity supersolution of (\ref{lemma4}). Now
for  $\lambda>0$ fixed, we define the nondecreasing
continuous real function
\[
\psi_{\lambda}(t)= \left\{\begin{array}{ll} \mu , & t\leq \mu ,
\\ \lambda^2 (t-\mu )+\mu , & \mu\leq t
\leq \mu+\frac{1}{\lambda} , \\ \lambda+\mu  & t\geq
\mu+\frac{1}{\lambda}.
\end{array} \right.
\]
The  function $\psi_{\lambda}\circ V$ is a viscosity supersolution
of equation (\ref{lemma4}) in $\mathcal{O}$ for every $\lambda$.
To prove this fact, we choose a sequence $\psi_n$ of strictly
increasing, smooth real maps that converge uniformly on compact
sets to $\psi_{\lambda}$. Then, for every $n$, the map
$\psi_n\circ V$ is a viscosity supersolution of equation
(\ref{lemma4}) in $\mathcal{O}$ by Lemma \ref{lemma2}. This
permits to conclude, by the stability of viscosity supersolutions
with respect to uniform convergence. Next we observe that the net
$\psi_{\lambda}\circ V$ is increasing and converges as
$\lambda\rightarrow +\infty$ to $V_\mu$. Viscosity supersolutions
are stable with respect to the pointwise increasing convergence
(see \cite{bcd}). Therefore the indicator function $V_\mu$ of
$K_\mu$ is a viscosity supersolution of equation (\ref{lemma4})
and then $K_\mu$ is viable for $\mu\leq\overline{\mu}$.

Conversely, we assume now that $K_\mu$ is viable for every $\mu\leq\overline{\mu}$.
Moreover, we observe that
\[\overline{V}(x)=\inf_{\mu\leq \overline{\mu}}V_\mu(x)\wedge \overline{\mu}=
\inf\{\mu\leq\overline{\mu}\ |\ V(x)\leq\mu\}\wedge\overline{\mu}.\]
So, by the stability properties of viscosity supersolutions, if
$V_\mu$  solves (\ref{lemma4}) for every $\mu\leq\overline{\mu}$, then $\overline{V}$
solves (\ref{lemma4}) too.
\end{proof}
This lemma provides the main tool to prove the direct Lyapunov
theorems (see \cite{bc2}, \cite{thesis}).
\begin{teo}[Direct Lyapunov theorem]Assume
(\ref{condition1}), (\ref{convex}).  If the system admits a
Lyapunov function in an open set $\mathcal{O}$ containing the
equilibrium then

\noindent (i) the system $(CSDE)$ is almost surely stabilizable at
the origin;

\noindent (ii) if, in addition, the domain ${\cal O}$ can be
chosen as  $\R^N$, the system is also a.s.  Lagrange stabilizable
and for all $x\in\R^N$ there exists $\oa_.\in {\cal A}_x$ such
that the corresponding trajectory $\ox_.$ satisfies
\begin{equation}
\label{stimaglob} \vert \ox_t \vert \leq
\gamma_1^{-1}(\gamma_2(|x|))\quad \forall\, t\geq 0\quad{a.s.}
\end{equation}
with $\gamma_1, \gamma_2\in {\cal K}_{\infty}$.\\
Assume moreover that $V$ is a strict Lyapunov function then

\noindent (i) the system $(CSDE)$ is a.s. locally asymptotically
stabilizable at the origin;

\noindent (ii) if, in addition, the domain ${\cal O}$ of $V$ is
all $\R^N$, the system is a.s. globally asymptotically
stabilizable.
\end{teo}
\section{Converse Lyapunov theorems for a.s. stabilizability}
In this section we prove the main results in the article: we
assume that the system (CSDE) satisfies an a.s. stabilizability
property and construct an appropriate Lyapunov function.
\begin{teo}[a.s. stabilizability]\label{asteoconverse1}
Assume (\ref{condition1}), (\ref{convex}). Then

\noindent (i) if the system $(CSDE)$ is almost surely stabilizable
at the origin in the ball $B_K$, the function \[V(x)=
\left[\inf_{\alpha\in \mathcal{A}_x}\sup_{t\geq 0}\esup \
|X_t^\alpha|\right]\wedge K\] is a Lyapunov function for the
system in $B_K$;

\noindent (ii) if the system is also a.s. Lagrange stabilizable
then the function \[V(x)= \inf_{\alpha\in
\mathcal{A}_x}\sup_{t\geq 0}\esup \ |X_t^\alpha|\]is a global
Lyapunov function for the system.
\end{teo}
\begin{proof}
We start proving (i). We start constructing the larger viable set
containing the origin in the ball $\overline{B_K}$. We construct a
nonnegative, uniformly continuous, radial function $c_K$ such that
$c_K(|s|)=0$ if $|s|\leq K$ and $0<c_K(|s|)\leq |s|$ for $|s|> K$
and for $\lambda<0$ we consider the value function
\[W_K(x)=\inf_\alpha\mathbf{E}_x\int_0^{+\infty}\!\!c_K(|X_s^\alpha|)
e^{-\lambda s}ds.\] It is well known (\cite{pll1} and \cite{fso})
that $W_K$ is a continuous viscosity supersolution of the
Hamilton-Jacobi-Bellman equation in $\R^N$
\[\max_{\alpha\in A}\left\{-DW(x)\cdot
f(x,\alpha)-trace\left[a(x,\alpha) D^2W(x)\right] \right\}+\lambda
W(x)\geq c_K(|x|) \] The function $W_K$ is nonnegative in $\R^N$
and   $W_K(0)=0$ since the origin is a controlled equilibrium as
remarked in (\ref{equilibrium}). We consider the propagation set
of the minimum value $0$:
\[Prop(0,W_K)=\{x\in\R^N\ |\ W_K(x)=0\}=\{x\ |\ \exists \oa \ :\
\ox_t\in \overline{B_K}\ a.s.\forall t\}\subseteq
\overline{B_K}.\] This set is clearly closed and it can be proved
that it is also viable. This follows immediately from Theorem 4.6
 in Bardi,  Da Lio \cite{bdl3} or can be checked directly.

The candidate Lyapunov function is defined in the set
$\overline{\mathcal{O}}=\overline{B_K}$ as
\[V(x)=\left\{\begin{array}{ll} \inf_{\alpha\in
\mathcal{A}_x}\sup_{t\geq 0}\esup\ |X_t^\alpha|& x\in
Prop(0,W_K)\\  K & x\in \overline{B_K}\setminus
Prop(0,W_K).\end{array}\right.\] First of all we observe that the
function $V$ is well defined: by definition of a.s.
stabilizability, there exists a function $\gamma\in\mathcal{K}$
such that $V(x)\leq \gamma(|x|)$. From this we get also  that  $V$
is continuous at the origin. We observe also that
$Prop(0,W_K)\supseteq B(0,\gamma^{-1}(K))$. Moreover  $V$ is
positive definite. Indeed if $V(x)=0$ then for every $\eps>0$
there exists $\alpha_\eps$ such that the corresponding trajectory
satisfies  $|X_t|\leq \eps$ for all $t\geq 0$ almost surely: so
$\inf_{\alpha_.\in {\cal A}_x} \mathbf{E}_0\int_0^{+\infty} |X_t|
e^{-\lambda t}dt=0 $ for any fixed $\lambda>0$. By Theorem
\ref{optimalcontrol},  the $\inf$ is attained, and the minimizing
control produces a trajectory satisfying a.s. $|X_t|=0$ for all
$t\geq 0$.

To prove the semicontinuity and the differential inequality, we
provide another characterization of  $V$. We prove that it
coincides with
\begin{equation}\label{wo}
w(x):=\inf\{r\ |\ \exists\oa\in \mathcal{A}_x \ |\ox_t|\leq r\ \
\forall t\geq 0 \ \ a.s.\}\wedge K.\end{equation} First of all
with the usual argument based on Theorem  \ref{optimalcontrol}, we
get that this infimum is actually a minimum. Then we fix $K\geq
k>w(x)$: by definition there exists $\oa\in\mathcal{A}_x$ such
that $|\ox_t|\leq k-\eps$ for all $t\geq 0$ a.s. This means that
$\sup_{t\geq 0}\mbox{ess}\sup_{\omega\in \Omega}|\ox_t|<k$ from
which we deduce $w(x)\geq V(x)$. The converse is similar.

Therefore, for every  $0\leq\mu < K$, the  sublevel set
$\{x\in\mathcal{O}\ |\ V(x)\leq \mu\}$ coincides with $\{x\in\R^N\
|\ \exists\oa\in \mathcal{A}_x \ |\ox_t|\leq \mu\ \ \forall t\geq
0 \ \ a.s.\}$; in particular this gives that the function $V$ is
proper. As at the beginning of the proof, we can characterize
these sets as the propagation sets of the minima of viscosity
supersolutions of suitable Hamilton-Jacobi-Bellman equations.
Then, by the results in \cite{bdl3}, these sets are closed in
$\R^N$ and viable with respect to (CSDE). So the function $V$ is
LSC. Moreover, by Lemma \ref{lemma3}, it satisfies the
infinitesimal decrease condition (\ref{cseq1}).

To prove (ii),  we observe that for every $K>0$ we can repeat the
previous construction, since the estimate (\ref{cappa}) holds in
the whole space with $\gamma\in\mathcal{K}_{\infty}$. So we get an
increasing sequence of Lyapunov functions $V_K$: for every $K>0$
we construct as before the function $V_K$ in the ball $B_K$ and
extend it to the whole space in the obvious way. Hence the global
Lyapunov function for the system is
\[V(x)=\lim_{K\to+\infty}V_K(x)=\sup_{K>0}V_K(x)=\inf_{\alpha\in \mathcal{A}_x}\sup_{t\geq
0}\mbox{ess}\sup_{\omega\in \Omega} |X_t^\alpha|.\] It is
immediate to check that it satisfies the condition (i),(ii),(iii)
in the Definition \ref{cliap}. Moreover the proof of the fact that
$V$ satisfies the differential condition (iv) relies on standard
stability properties of viscosity supersolutions  with respect to
the pointwise increasing convergence (\cite{bcd}).
\end{proof}
Now we prove the converse Lyapunov theorem for asymptotic
stability in a bounded set.
\begin{teo}[a.s. asymptotic stabilizability] \label{asteoconverse2}
Assume (\ref{condition1}), (\ref{convex}). If the system $(CSDE)$
is a.s.  asymptotically stabilizable at the origin in the ball
$B_K$, then there exist an open set $\mathcal{O}$ containing the
origin and a Lipschitz continuous positive definite function
$l:\mathcal{O}\to \R$ such that
\[V(x)=\inf_{\alpha\in\mathcal{A}_x}\esup
\int_0^{+\infty}l(X_t^\alpha)dt\] is a Lyapunov function for the
system in $\mathcal{O}$.
\end{teo}
\begin{proof}
The estimate (\ref{decayest.})  permits to construct a positive
definite Lipschitz continuous function $l:\R^N\rightarrow \R$ as
follows. $t_x(U)$ denotes the random time spent in the set $U$ by
the trajectory $\ox$. Using the properties of the function $\beta$
we get that, for every $r\leq K$ and $x\in B_K$, the time
$t_x(\R^N\setminus B_r)(\omega)$ spent by the trajectory $\ox$
outside $B_r$ is almost surely bounded and moreover it satisfies
\[\sup_{x\in B_K}\esup
t_x(\R^N\setminus B_r)(\omega)<+\infty.\] We consider now a
decreasing sequence of positive numbers $r_i$ such that $r_0<K$
and $\lim_{i\rightarrow+\infty}r_i=0$. We define
\[T_i=\sup_{x\in B_K}\esup t_x(\R^N\setminus
B_{r_i})(\omega).\] The sequence of positive numbers $T_i$ is
increasing as $i\rightarrow +\infty$: we can choose a decreasing
sequence of positive numbers $l_i$ such that
$\sum_{i=0}^{+\infty}l_iT_i=M<+\infty$. The function
$l:\R^N\rightarrow \R$ is therefore defined as  a radial Lipschitz
continuous, positive definite, nondecreasing function which
satisfies $l(0)=0$, $l(|x|)=l_{i+1}$ for $|x|=r_i$ and
$l(|x|)=l_1$ for every $|x|\geq r_0$.

The candidate Lyapunov function is
\[V(x)=\inf_{\alpha\in\mathcal{A}_x}\esup
\int_0^{+\infty}l(|X_t^\alpha|)dt\qquad \text{ for }x\in B_K.\]
The rest of the proof will be devoted to show that this function
satisfies the properties of Definition \ref{cliap} and then is a
strict Lyapunov function. First of all  $V$ is well defined:
\[V(x)=\inf_{\alpha\in\mathcal{A}_x}\esup
\int_0^{+\infty}l(|X_t^\alpha|)dt\leq
\esup\int_0^{+\infty}l(|\ox_t|)dt\leq
\sum_{i=0}^{+\infty}l_iT_i=M.\]  By definition $V(x)\geq 0$ for
every $x$ and $V(0)=0$. We assume now that for some $x\neq 0$
$V(x)=0$: this means that for every $\eps>0$ there exists
$\alpha_\eps\in\mathcal{A}_x$ such that $\int_0^{+
\infty}l(|X_t^{\alpha_\eps}|)dt\leq\eps$ almost surely. Then
$\inf_\alpha\mathbf{E}_x\int_0^{+\infty}l(|X_t^\alpha|)dt=0$: from
this, by the usual argument based on Theorem \ref{optimalcontrol},
we deduce that $x=0$. We show now that the function $V$ is
continuous at the origin. Recalling the definition of
$t_x(\R^N\setminus B_r)$ and the a.s. asymptotic stabilizability,
we get that $t_x(\R^N\setminus B_r)=0$ almost surely for initial
data $x$ such that $\beta(|x|,0)\leq r$. Therefore, by the
continuity at $r=0$ of the function $\beta(r,0)$, for every
$\eps>0$ there exists $\delta>0$ such that, for $|x|\leq \delta$,
$\beta(|x|,0)\leq \eps$ and then $t_x(\R^N\setminus B_k)=0$ almost
surely for every $k\geq\eps$. So we get $V(x)\leq
\sum_{i=i(\eps)}^{+\infty}l_i T_i$ where $i(\eps)$ is the minimum
index for which  $r_{i(\eps)}\leq \eps$. Since the sum $\sum_i
l_iT_i$ converges and $r_i\rightarrow 0$ as $i\rightarrow
+\infty$, for every $\theta>0$, we can choose $\eps>0$ such that
$\sum_{i=i(\eps)}^{+\infty}l_i T_i\leq\theta$: this gives the
continuity at the equilibrium. To conclude the proof we have to
provide another equivalent definition of $V$. We consider  the new
system   in $\R^{N+1}$ \[ (CSDE2)\left\{\begin{array}{l}
d(X_t,Y_t)=\overline{f}(X_t,Y_t,\alpha_t)dt+
\overline{\sigma}(X_t,Y_t,\alpha_t)d(B_t,0),\;\; t>0,\\
(X_0,Y_0)=(x,y).
\end{array}
\right.
\]
where $\overline{f}(x,y,\alpha)=(f(x,\alpha),l(x))$ and
$\overline{\sigma}(x,y,\alpha)= (\sigma(x,\alpha),0)$. It
satisfies  conditions (\ref{condition1}) and  (\ref{convex}) and
has $(0,0)$ as a controlled equilibrium. We introduce now the
following function
\[W(x,y)=\inf \{r\ |\
\exists \alpha\in\mathcal{A}_x \ |Y_t^\alpha|\leq r \ \forall
t\geq 0\text{ a.s.}\}=
\]\[=\inf\{r\ |\ \exists \alpha\in\mathcal{A}_x\
\left|y+\int_0^{+\infty}l(|X_t^\alpha|)dt\right|\leq r\text{
a.s.}\}.\]  Using Theorem \ref{optimalcontrol}, it can be proved
easily that this infimum is actually a minimum. So the sublevel
sets of $W$ are
\[\left\{(x,y)\ |\ W(x,y)\leq \mu\right\}=\left\{(x,y)\ |\ \exists
\alpha\in\mathcal{A}_x\
\left|y+\int_0^{+\infty}l(|X_t^\alpha|)dt\right|\leq \mu\text{
a.s.}\right\}.\]Repeating the argument in the proof of Theorem
\ref{asteoconverse1}, we get that these sets are closed and viable
with respect to (CSDE2). So the function $W(x,y)$ is LSC and
satisfies, by Lemma \ref{lemma3}, in viscosity sense
\begin{equation}\label{hbsinto}
\sup_{\sigma(x,\alpha)\cdot D_xW=0}\left\{-D_xW\cdot
f(x,\alpha)-trace[a(x,\alpha)D^2_{xx}W]\right\}-l(|x|)D_yW \geq 0.
\end{equation}
Now we show that the candidate Lyapunov function $V$ coincides
with the function $W$ on the set $B_K\times\{y=0\}$. Assume that
$V(x)\leq r$. For every $\eps>0$, there exists
$\alpha_\eps\in\mathcal{A}_x$ such that almost surely
$\int_0^{+\infty}l(|X_t^{\alpha_\eps}|)dt\leq r+\eps$. Therefore
$W(x,0)\leq r+\eps$ and so we conclude by the arbitrariness of
$\eps$. The proof of the opposite inequality $V(x)\leq W(x,0)$ is
similar. From this characterization of the function $V$ we deduce
immediately that it is LSC and bounded in $B_K$. It remains to
check the differential condition (\ref{cseq2}). We fix $x\neq 0$
in $B_K$ and consider $(p,Y)\in\mathcal{J}^{2,-}V(x)$: by
definition, for every $x'\rightarrow x$, we get $W(x',0)\geq
W(x,0)+p\cdot(x'-x)+1/2(x'-x)Y(x'-x)+o(|x-x'|)$. Using the
definition, it is immediate to check $W(x',0)\leq W(x',y)-y$. This
implies  that if $(p,Y)\in\mathcal{J}^{2,-}V(x)$ then
$(p,Y,1)\in\mathcal{J}^{2,1,-}W(x,0)$. So the differential
condition (\ref{cseq2}) comes from (\ref{hbsinto}).
\end{proof}
\section{Extensions}
We can extend the results to the characterization of
stabilizability of a closed set $M\subseteq\R^N$. We denote with
$d(x,M)$ the distance between $x\in \R^N$ and $M$.
 \begin{defi}[a.s.  stabilizability  at $M$]\label{def:mstab}
The system $(CSDE)$ is \emph{almost surely
(open loop) stabilizable} at $M$ if % it is viable for the dynamical
there exists $\gamma\in\mathcal{K}$ such
that, for every $x$ in a neighborhood of $M$,
there is an admissible control function
 $\overline{\alpha}_{\cdot}\in{\cal A}_x$
 %with the property that the corresponding
 whose trajectory
 $\overline{X}_{\cdot}%^{\overline{\alpha}_{\cdot}},
$
%starting at $x$
verifies
$$
d(\overline{X}_t,M)\leq \gamma(d(x,M)) \quad
\forall \, t\geq 0 \quad
%\text{a.s.}
\text{almost surely.}
$$
\end{defi}
\begin{rem}\upshape
From this definition, using Theorem \ref{optimalcontrol}, we can
deduce that if $M$ is a.s.  stabilizable, then it is viable for
$(CSDE)$.
\end{rem}
We adapt the definition of control Lyapunov function to the case the
equilibrium is a set $M$:
 \begin{defi}[control Lyapunov functions at $M$]
\label{cmliap} Let ${\cal O}$ be an open neighborhood of the
closed set $M$. A function $V:{\cal O} \rightarrow[0,+\infty)$ is
a control Lyapunov function at $M$ for $(CSDE)$ if

\noindent (i) $V$  is lower semicontinuous;

\noindent(ii)
there exist $\gamma_1, \gamma_2 \in\mathcal{K}$ such that
$
\gamma_2(d(x,M))\leq V(x)\leq\gamma_1 (d(x,M))
%\quad \forall x\in\R^N
$
for all $x\in{\cal O}$;

\noindent (iii) for all $x\in{\cal O}\setminus M$ and
$(p,Y)\in{\cal J}^{2,-}V(x)$ there exists $\overline{\alpha}\in A$
such that condition (\ref{cseq2}) holds.
\end{defi}
\begin{teo}
Assume (\ref{condition1}), (\ref{convex}). Then the system
$(CSDE)$ is almost surely stabilizable at $M$ if and only if
 there exists a
control Lyapunov function $V$ at $M$.
\end{teo}

\end{document}